\documentclass[a4paper,reqno,10pt]{amsart}
\usepackage{amssymb}
\usepackage{cite}
\usepackage{tikz}

\usepackage{amsfonts}
\usepackage{bigints}
\usepackage{amsmath}
\usepackage{amssymb}
\usepackage{systeme}
\usepackage{mathrsfs}
\usepackage[colorlinks=true]{hyperref}
\numberwithin{equation}{section}
\newtheorem{theorem}{Theorem}[section]
\newtheorem{prop}[theorem]{Proposition}
\newtheorem{lm}[theorem]{Lemma}
\newtheorem{cor}[theorem]{Corollary}

\newtheorem{rem}{Remarks}
\newcommand{\R}{\mathbb{R}}

\def\dive{\mathrm{div}}
\newenvironment{preuve}{{\noindent {\bf Proof. }}}{\hfill {\rule{2.5mm}{2.5mm}}}
\newenvironment{preuve1}{{\noindent {\bf Proof of Theorem \ref{theorem1}. }}}{\hfill {\rule{2.5mm}{2.5mm}}}

\makeatletter
\let\@msm@th@eqref\eqref
\renewcommand{\eqref}[1]{%
	\begingroup
	\leavevmode
	\color{red}%
	\hypersetup{linkbordercolor=[named]{red}}%
	\@msm@th@eqref{#1}%
	\endgroup
}
\makeatother

\author[M.~Amara]{Mustapha Amara}
\address{Department of Mathematics, Faculty of Science of Gab\`es, Research Laboratory Mathematics and Applications LR17ES11; Tunisia}
\email{\sl Mostafa.Amara@fsg.u-gabes.tn}

\title[On the boundedness of the global solution of $(AQG)$ equations in Sobolev space]
{On the boundedness of the global solution of anisotropic quasi-geostrophic equations in Sobolev space}

\begin{document}
	\begin{abstract}
	In this paper, we show that the global solution of the surface anisotropic two-dimensional quasi-geostrophic equation with fractional horizontal dissipation and vertical thermal diffusion established by the author in \cite{YZ} is bounded in Sobolev spaces uniformly with respect to time.
	\end{abstract}
	
	

	\maketitle
	\tableofcontents

	
	\section{\bf Introduction}
	In this paper we deal with the following surface quasi-geostrophic  equation with fractional horizontal dissipation and fractional vertical thermal diffusion:
	\begin{equation}\label{AQG}\tag{AQG}
		\begin{cases}
			\partial_t\theta+ u_\theta.\nabla\theta +\mu|\partial_1|^{2\alpha}\theta+\nu |\partial_2|^{2\beta}\theta=0,&x\in\R^2, \ t>0\\
			\theta(x,0)=\theta^0(x),
		\end{cases}
	\end{equation}
where $\alpha\in(0,1)$, $\beta\in (0,1)$, $\mu>0$ and $\nu>0$ are real number. Here, we denote by $|\partial_1|$ and  $|\partial_2|$ the operators given by 
$$\mathcal{F}(|\partial_1|^{2\alpha}f)(\xi)=|\xi_1|^{2\alpha}\mathcal{F}(f)(\xi),\quad\mathcal{F}(|\partial_2|^{2\beta}f)(\xi)=|\xi_2|^{2\beta}\mathcal{F}(f)(\xi),\quad\forall \xi=(\xi_1,\xi_2)\in \R^2,$$
where $\mathcal{F}(f)$ represents the Fourier transformation of $f$. The variable $\theta$ represents the potential temperature and the velocity $u_\theta=(u_1,u_2)$ is determined by $\theta$ via the formula
	\begin{align*}
		u_\theta=\mathcal{R}^\perp\theta=\left( -\mathcal{R}_2\theta,\mathcal{R}_1\theta\right)=\left(-\partial_{2}(-\Delta)^{-\frac{1}{2}}\theta,\partial_{1} (-\Delta)^{-\frac{1}{2}}\theta\right),
	\end{align*}
	where $\mathcal{R}_1,\mathcal{R}_2$ are the standard 2D Riesz transforms. Clearly, the velocity $u$ is divergence free,
	namely $$\dive(u_\theta):=\partial_{1}u_1+\partial_{2}u_2=0.$$
	
	The system \eqref{AQG} is deeply related, in the case when $\alpha=\beta$ and $\mu=\nu$,  to the classical dissipation $(QG)$ equation, with its from as follows
		\begin{equation}\label{QG}\tag{QG}
		\begin{cases}
			\partial_t\theta+ u_\theta.\nabla\theta +\mu(-\Delta)^{2\alpha}\theta=0,&x\in\R^2, \ t>0\\
			\theta(x,0)=\theta^0(x).
		\end{cases}
	\end{equation}
	
	The equations \eqref{AQG} and \eqref{QG} are special cases of the general quasi-geostrophic approximations for atmospheric and oceanic fluid flow with small
	Rossby and Ekman numbers. The first mathematical studies of this equation was carried out in 1994s by Constantin, Majda and Tabak. For more
	details and mathematical and physical explanations of this model we can consult \cite{DC,CP,JP,CP1}.\\
	
The inviscid quasi-geostrophic equation (i.e., \eqref{QG} with $\mu=0$) shares many properties parallel to those of 3D Euler equations such as the vortex stretch mechanism and thus serves as a lower dimensional model of 3D Euler equations.\\
	
	The first studies of the system \eqref{AQG} is in \cite{YZ} by Ye, who shows that this equation admits a unique global solution $\theta$ in the space $C(\R^+,H^s(\R^2))$, $s\geq 2$,  such that
	\begin{equation}
		|\partial_{1}|^\alpha\theta,|\partial_{2}|^\beta\theta\in L^2_{loc}(\R^+,H^s(\R^2)),
	\end{equation} when $\alpha,\beta\in (0,1)$ satisfies
	\begin{equation}\label{1.1}
		\beta>\begin{cases}
			\frac{1}{2\alpha+1},&0<\alpha\leq \frac{1}{2}\\ \\
			
			\frac{1-\alpha}{2\alpha}&\frac{1}{2}<\alpha<1.
		\end{cases}
	\end{equation}

	We also refer to our result in \cite{MJ}, when we established the global regularity in $C(\R^+,H^s(\R^2))$, when $\alpha,\beta\in (1/2,1)$ and $s\in (\max\{2-2\alpha,2-2\beta\},2)$.\\
	
In this paper, we will show that this global solution is in $L^\infty(\R^+,H^s(\R^2))$ and
\begin{equation}
	|\partial_{1}|^\alpha\theta,|\partial_{2}|^\beta\theta\in L^2(\R^+,H^s(\R^2)).
\end{equation}

	For the sake of simplicity, we will set $\mu=\nu=1$ throughout the paper, $C$ denoted all constants that is a generic constant depending only on the quantities specified in the context and $C(\theta^0)$ represent all constants depending on the norm of the initial condition $\theta^0$.
	
	\section{\bf Main theorems}
We explain all the details later in the paper, but let us state here the first main result:
\begin{theorem}\label{theorem1}
	Let $\theta^0\in H^s(\R^2)$ for $s\geq 2$. If $\alpha,\beta\in (0,1)$ satisfy \eqref{1.1}. Then the system \eqref{AQG} admits a unique global solution $\theta$ such that
	$$\theta\in C_b(\R^+,H^s(\R^2)),\quad |\partial_{1}|^\alpha\theta,|\partial_{2}|^\beta\theta\in L^2(\R^+,H^s(\R^2)).$$
	Moreover, for any $t\geq 0$
 \begin{equation}
	\|\theta(t)\|^2_{H^s}+\int_0^t||\partial_{1}|^{\alpha}\theta\|^2_{H^s}d\tau+\int_{0}^t\||\partial_2|^{\beta}\theta\|^2_{H^s}d\tau\leq C(\theta^0).
\end{equation}
	\end{theorem}
 We outline the main ideas in the proof of this Theorem which indicated by the author in \cite{YZ} as follows:
 \begin{theorem}[see \cite{YZ}]\label{thm2}
 	Let $\theta^0\in H^s(\R^2)$ for $s\geq 2$. If $\alpha,\beta\in (0,1)$ satisfy \eqref{1.1}.
 	Then the system \eqref{AQG} admits a unique global solution $\theta$ such that  for any $T>0$ we have 
 	\begin{equation}
 		\theta\in C([0,T],H^s(\R^2)),\ |\partial_{1}|^{\alpha}\theta,|\partial_{2}|^{\beta}\theta\in L^2([0,T],H^s(\R^2)).
 	\end{equation}
 Moreover, for any $t\geq 0$
  \begin{align}
 	\label{Lp}	 &\hskip6cm	\|\theta(t)\|_{L^p}\leq \|\theta^0\|_{L^p},\hskip4cm\forall p\in [2,+\infty],\\
 	&	\label{L2}\|\theta(t)\|_{L^2}+\int_{0}^{t}\||\partial_{1}|^\alpha\theta\|_{L^2}^2d\tau+ \int_{0}^{t}\||\partial_{2}|^\beta\theta\|_{L^2}^2d\tau\leq \|\theta^0\|_{L^2}^2,\\
 	\label{H1}	 &	\|\theta(t)\|^2_{\dot{H}^1}+\int_{0}^{t}\||\partial_{1}|^{\alpha}\theta\|^2_{\dot{H}^1}d\tau+\int_{0}^{t}\||\partial_2|^{\beta}\theta\|^2_{\dot{H}^1}d\tau\leq\|\theta^0\|^2_{\dot{H}^1}+ C\int_{0}^{t} \left(1+\|u\|_{L^\infty}^{\rho}\right)\|\theta\|^2_{\dot{H}^1}d\tau,\\
 	\nonumber &	\|\theta(t)\|^2_{\dot{H}^2}+\int_{0}^{t}\||\partial_{1}|^{\alpha}\theta\|^2_{\dot{H}^2}d\tau+\int_{0}^{t}\||\partial_2|^{\beta}\theta\|^2_{\dot{H}^2}d\tau\\
 	\label{H2}	 &\hskip5cm\leq\|\theta^0\|_{\dot{H}^2}+C\int_{0}^{t} \left(1+ \|\theta\|_{\dot{H}^1}^{2}\right)\left(1+\||\partial_{1}|^\alpha\theta\|_{\dot{H}^1}^2+\||\partial_{2}|^\beta\theta\|_{\dot{H}^1}^2\right)\|\theta\|_{\dot{H}^2}^2d\tau,
 \end{align}
 and 
 \begin{equation}\label{Hs}
 	\|\theta(t)\|^2_{H^s}+\int_0^t||\partial_{1}|^{\alpha}\theta\|^2_{H^s}d\tau+\int_{0}^t\||\partial_2|^{\beta}\theta\|^2_{H^s}d\tau\leq \|\theta^0\|^2_{H^s}+ C\int_{0}^t\left(1+\|\nabla u\|_{L^\infty}+\|\nabla \theta\|_{L^\infty}\right) \| \theta\|_{H^s}^2d\tau,
 \end{equation}
 where 
 \begin{equation}\label{rho}
 	\rho=\begin{cases}
 		\frac{2\beta}{(2\alpha+1)\beta-1},&\mbox{ if }\beta>\frac{1}{2\alpha+1},\ \alpha\leq \frac{1}{2},\\
 		\\
 		\max\left\{\frac{2\alpha}{2\alpha-1},\frac{2\alpha}{(2\beta+1)\alpha-1}\right\},& \mbox{if } \beta>\frac{1-\alpha}{2\alpha},\ \alpha> \frac{1}{2}.
 	\end{cases}
 \end{equation}
 \end{theorem}
\begin{rem}
	Our second result is a corollary of the first main result and the proved result in \cite{MJ}, which show that for any initial condition $\theta^0$ in $H^s(\R^2)$, $s\in (\max\{2-2\alpha,2-2\beta\},2)$, $\alpha,\beta$ in $(1/2,1)$, we have a global solution $\theta$ of the system \eqref{AQG} satisfy
	\begin{equation}
		\theta\in C(\R^+,H^s(\R^2)),\ |\partial_{1}|^{\alpha}\theta,|\partial_{2}|^{\beta}\theta\in L^2_{loc}(\R^+,H^s(\R^2)).
	\end{equation}
Moreover, there is $T_0>0$, such that
$$\theta\in C((0,T_0),H^2(\R^2)),$$
which implies that $\theta\in C((0,+\infty),H^2(\R^2)).$
\end{rem}
Now, we ready to state our second result:
\begin{cor}\label{Corollary2.5}
	Let	$\alpha,\beta\in (1/2,1)$ and  $\theta^0\in H^s(\R^2)$, $s\in (\max\{2-2\alpha,2-2\beta\},2) $ be the global solution of \eqref{AQG}. Then, the system \eqref{AQG} admits a unique global solution $\theta$ such that
	$$\theta\in C_b(\R^+,H^s(\R^2)),\quad |\partial_{1}|^\alpha\theta,|\partial_{2}|^\beta\theta\in L^2(\R^+,H^s(\R^2)).$$
\end{cor}
\section{\bf Notations and Preliminary Results}
	In this short section, we collect some notations and definitions that will be used later, and we	give some technical lemmas.
	\subsection{\bf Notations}
	\begin{itemize}
		\item[$\bullet$] The Fourier transformation in $\R^2$		
		\begin{equation}
			\mathcal{F}(f)(\xi)=\widehat{f}(\xi)=\int_{\R^2}e^{-ix.\xi}f(x)dx,\quad \xi\in \R^2.
		\end{equation}
		The inverse Fourier formula is
		\begin{equation}
			\mathcal{F}^{-1}(f)(x)=(2\pi)^{-2}\widehat{f}(\xi)=\int_{\R^2}e^{i\xi.x}f(\xi)d\xi,\quad x\in \R^2.
		\end{equation}
		\item[$\bullet$] The convolution product of a suitable pair of function $f$ and $g$ on $\R^2$ 	is given by
		\begin{equation}
			f\ast g(x)=\int_{\R^2} f(x-y)g(y)dy
		\end{equation}
		\item[$\bullet$] If $f=(f_1,f_2)$ and $g=(g_1,g_2)$ are two vector fields, we set
		$$f\otimes g:= (g_1f,g_2f),$$
		and 
		$$\dive (f\otimes g):= (\dive(g_1f),\dive(g_2f)).$$
		\item[$\bullet$] For $s\in \R$:
		\begin{itemize}
			\item[$\ast$] The Sobolev space: 
			\begin{align*}
				&H^{s}(\R^2):=\left\{f\in \mathcal{S}'(\R^2); (1+|\xi|^2)^{s/2}\widehat{f}\in L^2(\R^2)\right\},
			\end{align*} denotes the usual inhomogeneous Sobolev space on $\R^2$, with the norm 
			\begin{align*}
				&	\|f\|_{H^s}=\left(\int_{\R^2}(1+|\xi|^2)^s|\widehat{f}(\xi)|^2d\xi\right)^{\frac{1}{2}}
			\end{align*}
		\item[$\ast$] The Sobolev space: 
		\begin{align*}
			 \dot{H}^{s}(\R^2):=\left\{f\in \mathcal{S}'(\R^2);\widehat{f}\in L^1(\R^2)\mbox{ and }|\xi|^s\widehat{f}\in L^2(\R^2)\right\},
		\end{align*} denotes the usual  homogeneous Sobolev space on $\R^2$, with the norm 
		\begin{align*}
			\|f\|_{\dot{H}^s}=\left(\int_{\R^2}|\xi|^{2s}|\widehat{f}(\xi)|^2d\xi\right)^{\frac{1}{2}}.
		\end{align*}
		\end{itemize}
	\end{itemize}
\subsection{\bf Preliminary Results}
	We recall a fundamental lemma concerning the Sobolev spaces
	\begin{lm}\label{Lemma1}
		For $s,s_1,s_2\in\R$ and $z\in[0,1]$, the following anisotropic interpolation inequalities hold true for $i=1,2$:
		\begin{align}
		\label{ing1}	&\||\partial_{i}|^{zs_1+(1-z)s_2}f\|_{H^{s}}\leq \||\partial_{i}|^{s_1}f\|_{H^{s}}^z\||\partial_{i}|^{s_2}f\|_{H^{s}}^{1-z},\\
		\label{ing2}	&\||\partial_{i}|^{zs_1+(1-z)s_2}f\|_{\dot{H}^{s}}\leq \||\partial_{i}|^{s_1}f\|_{\dot{H}^{s_1}}^z\||\partial_{i}|^{s_2}f\|_{\dot{H}^{s_2}}^{1-z}.
		\end{align}
	\end{lm}
\begin{preuve}
It suffices to show \eqref{ing1} for $i = 1$ as $i = 2$ can be performed as the same manner: So we have 
	\begin{align*}
		\||\partial_{i}|^{zs_1+(1-z)s_2}f\|_{H^{s}}^2&=\int_{\R^2}(1+|\xi|^2)^s|\xi_i|^{2(zs_1+(1-z)s_2)}\left|\widehat{f}(\xi)\right|^2d\xi\\
		&=\int_{\R^2}\left((1+|\xi|^2)^s|\xi_i|^{2s_1}\left|\widehat{f}(\xi)\right|^2\right)^z\left((1+|\xi|^2)^s|\xi_i|^{2s_2}\left|\widehat{f}(\xi)\right|^2\right)^{1-z}d\xi=\|f_1\times f_2\|_{L^1},
	\end{align*}
where 
$$f_1(\xi)=\left((1+|\xi|^2)^s|\xi_i|^{2s_1}\left|\widehat{f}(\xi)\right|^2\right)^z,\quad f_2(\xi)=\left((1+|\xi|^2)^s|\xi_i|^{2s_2}\left|\widehat{f}(\xi)\right|^2\right)^{1-z}.$$
The fact that 
$$\frac{1}{\dfrac{1}{z}}+\frac{1}{\dfrac{1}{1-z}}=1.$$
Using the Hölder's inequality, we get
\begin{align*}
	\|f_1\times f_2\|_{L^1}&\leq \|f_1\|_{L^{\frac{1}{z}}}\|f_2\|_{L^{\frac{1}{1-z}}}\\
	&\leq \left(\int_{\R^2}(1+|\xi|^2)^s|\xi_i|^{2s_1}\left|\widehat{f}(\xi)\right|^2d\xi\right)^z\left(\int_{\R^2}(1+|\xi|^2)^s|\xi_i|^{2s_2}\left|\widehat{f}(\xi)\right|^2d\xi\right)^{1-z}\\
	&\leq \||\partial_{i}|^{s_1}f\|_{H^{s}}^z\||\partial_{i}|^{s_2}f\|_{H^{s}}^{1-z},
\end{align*}
which implies the result for case inhomogeneous. Following the proof of \eqref{ing1}, the estimate \eqref{ing2} immediately holds true. This completes the proof of the lemma.
\end{preuve}
	\begin{lm}\cite{JN}\label{Lemma2}
		For any $p\in (1,+\infty)$, there is a constant $C(p)>0$ such that
		\begin{equation}
			\|\mathcal{R}^\perp\theta\|_{L^p}\leq C(p) \|\theta\|_{L^p}.
		\end{equation}
	\end{lm}
	\begin{lm}\label{Lemma3}\cite{YZ}
Let $\sigma>1$, then, their exist a constant $C>0$, such that for any $f\in H^\sigma(\R^2)$
	\begin{equation*}
		\|\mathcal{R}^\perp f\|_{L^\infty}\leq C\left(1+\|f\|_{L^2}+\|f\|_{L^\infty}\ln\left(e+\||\nabla|^\sigma f\|_{L^2}\right)\right) .
	\end{equation*}
\end{lm}

We finish with the following elementary inequality
	\begin{lm}\label{Lemma6}
		Let $\alpha>0$, then, there is a constant $C(\alpha)>0$ such that
		\begin{equation}
			\ln(x)\leq C(\alpha)x^\alpha,\quad\forall x\geq 1.
		\end{equation}
	\end{lm}
\section{Proof of Theorem \ref{theorem1}}
Before we begin the proof, we require the following proposition to prove that the global solution of system \eqref{AQG} provided in the Theorem \ref{thm2} is uniformly bounded in $H^1$ and $H^2$:
\begin{prop}\label{Prop1}
	If $\theta^0$, $\alpha$ and $\beta$ satisfies the assumptions stated in Theorem \ref{theorem1} and let $\theta$ be the corresponding global solution, then, for any $t\geq 0$,
	 \begin{equation}
		\|\theta(t)\|^2_{\dot{H}^1}+\int_0^t||\partial_{1}|^{\alpha}\theta\|^2_{\dot{H}^1}d\tau+\int_{0}^t\||\partial_2|^{\beta}\theta\|^2_{\dot{H}^1}d\tau\leq C(\theta^0).
	\end{equation}	
\end{prop}
\begin{preuve}
	We have the inequality \eqref{H1} given in Theorem \ref{thm2} given as following, for any $t\geq 0$,
	\begin{equation}
		\|\theta(t)\|^2_{\dot{H}^1}+\int_{0}^{t}\||\partial_{1}|^{\alpha}\theta\|^2_{\dot{H}^1}d\tau+\int_{0}^{t}\||\partial_2|^{\beta}\theta\|^2_{\dot{H}^1}d\tau\leq\|\theta^0\|^2_{\dot{H}^1}+ C\int_{0}^{t} \left(1+\|u\|_{L^\infty}^{\rho}\right)\|\theta\|^2_{\dot{H}^1}d\tau,
	\end{equation}
	where $\rho>1$ give in \eqref{rho}.\\
	
	In order to control $\|u_\theta\|_{L^\infty}=\|\mathcal{R}^\perp\theta\|_{L^\infty
	}$, we need the  logarithmic Sobolev interpolation inequality given in Lemma \ref{Lemma3}, so we obtain
	\begin{equation*}
		\|u_\theta\|_{L^\infty}\leq C\left(1+\|\theta\|_{L^2}+\|\theta\|_{L^\infty}\ln\left(e+||\nabla|^\sigma\theta\|_{L^2}\right)\right),
	\end{equation*}
	where $1<\sigma<1+\min\{\alpha,\beta\}$. Using \eqref{Lp} in the last inequality, we get
	\begin{equation*}
		\|u_\theta\|_{L^\infty}\leq C(\theta^0)\left(1+\ln\left(e+\||\nabla|^\sigma\theta\|_{L^2}\right)\right).
	\end{equation*}
	Moreover, by the Lemma \ref{Lemma6} we have  
	\begin{align*}
		\ln\left(e+\||\nabla|^\sigma\theta\|_{L^2}\right)&\leq C \left(1+\||\nabla|^\sigma\theta\|_{L^2}^{\frac{\sigma-1}{\sigma\rho}}\right).
	\end{align*}
	Finally, we get 
	\begin{align*}
		\|u_\theta\|_{L^\infty}^{\rho}&\leq C(\theta^0) \left(1+\||\nabla|^\sigma\theta\|_{L^2}^{\frac{\sigma-1}{\sigma}}\right)
	\end{align*}
	and 
	\begin{align*}
		\|\theta(t)\|^2_{\dot{H}^1}+\int_{0}^{t}\||\partial_{1}|^{\alpha}\theta\|^2_{\dot{H}^1}d\tau+\int_{0}^{t}\||\partial_2|^{\beta}\theta\|^2_{\dot{H}^1}d\tau&\leq\|\theta^0\|^2_{\dot{H}^1}+ C(\theta^0)\int_{0}^{t} \left(1+\||\nabla|^\sigma\theta\|_{L^2}^{\frac{\sigma-1}{\sigma}}\right)\|\theta\|^2_{\dot{H}^1}d\tau\\
		&\leq\|\theta^0\|^2_{\dot{H}^1}+ C(\theta^0)\int_{0}^{t}\|\theta\|^2_{\dot{H}^1}d\tau +C(\theta^0)\int_{0}^{t}\||\nabla|^\sigma\theta\|_{L^2}^{\frac{\sigma-1}{\sigma}}\|\theta\|^2_{\dot{H}^1}d\tau.
	\end{align*}
	Now, we need to control $\|\theta\|_{\dot{H}^1}$, so, we apply the the interpolation inequality. The fact that $0<1<\sigma$, then $1=z\times\sigma+(1-z)\times0$, $(z=\frac{1}{\sigma})$, and we have 
	\begin{align*}
		\|\theta\|_{\dot{H}^1}&\leq \|\theta\|_{L^2}^{1-z}\||\nabla|^\sigma\theta\|_{L^2}^{z}\\
		&\leq \|\theta^0\|_{L^2}^{\frac{\sigma-1}{\sigma}}\||\nabla|^\sigma\theta\|_{L^2}^{\frac{1}{\sigma}}.
	\end{align*}
	Therefore
	\begin{align*}
		C(\theta^0)\||\nabla|^\sigma\theta\|_{L^2}^{\frac{\sigma-1}{\sigma}}\|\theta\|^2_{\dot{H}^1}&\leq C(\theta^0) \||\nabla|^\sigma\theta\|_{L^2}^{\frac{\sigma-1}{\sigma}} \||\nabla|^\sigma\theta\|_{L^2}^{\frac{1}{\sigma}} \|\theta\|_{\dot{H}^1}\\
		&\leq C(\theta^0) \||\nabla|^\sigma\theta\|_{L^2} \|\theta\|_{\dot{H}^1}\\
		&\leq C(\theta^0) \|\theta\|^2_{\dot{H}^1}+ \||\nabla|^\sigma\theta\|_{L^2}^2,
	\end{align*}
	which imply
	\begin{align}\label{EST1}
		\|\theta(t)\|^2_{\dot{H}^1}+\int_{0}^{t}\||\partial_{1}|^{\alpha}\theta\|^2_{\dot{H}^1}d\tau+\int_{0}^{t}\||\partial_2|^{\beta}\theta\|^2_{\dot{H}^1}d\tau&\leq\|\theta^0\|^2_{\dot{H}^1}+ C(\theta^0)\int_{0}^{t}\|\theta\|^2_{\dot{H}^1}d\tau +\int_{0}^{t}\||\nabla|^\sigma\theta\|_{L^2}^2d\tau.
	\end{align}
	Since $\alpha<1<\sigma<1+\alpha$ and $\beta<1<\sigma<1+\beta$, then, using the interpolation inequality we get
	\begin{align}
		\nonumber	\||\nabla|^\sigma\theta\|^2_{L^2}&\leq 2 \||\partial_{1}|^\sigma\theta\|^2_{L^2}+2\||\partial_{2}|^\sigma\theta\|^2_{L^2}\\
		\label{EST2}	&\leq C\left(\||\partial_{1}|^\alpha\theta\|^2_{L^2}+\||\partial_{2}|^\beta\theta\|^2_{L^2}\right)+\frac{1}{4 }\left(\||\partial_{1}|^\alpha\theta\|^2_{\dot{H}^1}+\||\partial_{2}|^\beta\theta\|^2_{\dot{H}^1}\right); 
	\end{align} 
	and
	\begin{align}
		\nonumber C(\theta^0)	\|\theta\|^2_{\dot{H}^1}&\leq  C(\theta^0) \left(\||\partial_{1}|^\sigma\theta\|^2_{L^2}+ \||\partial_{2}|^\sigma\theta\|^2_{L^2}\right)\\
		\label{EST3}	&\leq C(\theta^0)\left(\||\partial_{1}|^\alpha\theta\|^2_{L^2}+\||\partial_{2}|^\beta\theta\|^2_{L^2}\right)+\frac{1}{4}\left(\||\partial_{1}|^\alpha\theta\|^2_{\dot{H}^1}+\||\partial_{2}|^\beta\theta\|^2_{\dot{H}^1}\right).
	\end{align} 
	Collect \eqref{EST2} and \eqref{EST3} with \eqref{EST1}, we get 
	\begin{align*}
		\|\theta(t)\|^2_{\dot{H}^1}+\int_{0}^{t}\||\partial_{1}|^{\alpha}\theta\|^2_{\dot{H}^1}d\tau+\int_{0}^{t}\||\partial_2|^{\beta}\theta\|^2_{\dot{H}^1}d\tau\leq\|\theta^0\|^2_{\dot{H}^1}&+ C(\theta^0)\int_{0}^{t}\left(\||\partial_{1}|^\alpha\theta\|^2_{L^2}+\||\partial_{2}|^\beta\theta\|^2_{L^2}\right)d\tau\\
		& +\frac{1}{2}\int_{0}^{t}\left(\||\partial_{1}|^\alpha\theta\|^2_{\dot{H}^1}+\||\partial_{2}|^\beta\theta\|^2_{\dot{H}^1}\right)d\tau.
	\end{align*}
	Finally we have 
	\begin{equation}
		\|\theta(t)\|^2_{\dot{H}^1}+\frac{1}{2}\int_{0}^{t}\||\partial_{1}|^{\alpha}\theta\|^2_{\dot{H}^1}d\tau+\frac{1}{2}\int_{0}^{t}\||\partial_2|^{\beta}\theta\|^2_{\dot{H}^1}d\tau\leq\|\theta^0\|^2_{\dot{H}^1}+ C(\theta^0) \|\theta^0\|^2_{L^2}.
	\end{equation}
	So we get that $$\theta\in C_b(\R^+,H^1(\R^2)),\ |\partial_{1}|^\alpha\theta,|\partial_{2}|^\beta\theta\in L^2(\R^+,H^1(\R^2)).$$
\end{preuve}
\begin{prop}
	If $\theta^0$, $\alpha$ and $\beta$ satisfies the assumptions stated in Theorem \ref{theorem1} and let $\theta$ be the corresponding global solution, then, for any $t\geq 0$,
	\begin{equation}
		\|\theta(t)\|^2_{\dot{H}^2}+\int_0^t||\partial_{1}|^{\alpha}\theta\|^2_{\dot{H}^2}d\tau+\int_{0}^t\||\partial_2|^{\beta}\theta\|^2_{\dot{H}^2}d\tau\leq C(\theta^0).
	\end{equation}	
\end{prop}\begin{preuve}
Using \eqref{H2} and Proposition \ref{Prop1}, we get for any $t\geq 0$,
\begin{align}
	\|\theta(t)\|^2_{\dot{H}^2}+\int_{0}^{t}\||\partial_{1}|^{\alpha}\theta\|^2_{\dot{H}^2}d\tau+\int_{0}^{t}\||\partial_2|^{\beta}\theta\|^2_{\dot{H}^2}d\tau\leq\|\theta^0\|_{\dot{H}^2}^2+C(\theta^0)\int_{0}^{t} \left(1+\||\partial_{1}|^\alpha\theta\|_{\dot{H}^1}^2+\||\partial_{2}|^\beta\theta\|_{\dot{H}^1}^2\right)\|\theta\|_{\dot{H}^2}^2d\tau.
\end{align}
Moreover, by the interpolation inequality for $\|\theta\|_{\dot{H}^2}$, where $\alpha<2<2+\alpha$ and $\beta<2<2+\beta$, we have
\begin{align*}
	C(\theta^0)	\|\theta\|^2_{\dot{H}^2}&\leq C(\theta^0)	\left(\||\partial_{1}|^2\theta\|^2_{L^2}+4\||\partial_{2}|^2\theta\|^2_{L^2}\right)\\
	&\leq C(\theta^0)\left(\||\partial_{1}|^\alpha\theta\|^2_{L^2}+\||\partial_{2}|^\beta\theta\|^2_{L^2}\right)+\frac{1}{2}\left(\||\partial_{1}|^\alpha\Delta\theta\|^2_{L^2}+\||\partial_{2}|^\beta\Delta\theta\|^2_{L^2}\right); 
\end{align*}
which implies
\begin{align}
	\|\theta(t)\|^2_{\dot{H}^2}+\int_{0}^{t}\||\partial_{1}|^{\alpha}\theta\|^2_{\dot{H}^2}d\tau+\int_{0}^{t}\||\partial_2|^{\beta}\theta\|^2_{\dot{H}^2}d\tau\leq C(\theta^0)\|\theta^0\|_{H^2}^2+C(\theta^0)\int_{0}^{t} \left(\||\partial_{1}|^\alpha\theta\|_{\dot{H}^1}^2+\||\partial_{2}|^\beta\theta\|_{\dot{H}^1}^2\right)\|\theta\|_{\dot{H}^2}^2d\tau.
\end{align}
Gronwall lemma implies that, for any $t\geq 0$ we have
$$ \|\theta(t)\|^2_{\dot{H}^2}\leq C(\theta^0)\|\theta^0\|_{H^2}^2\exp\left[ C(\theta^0) \int_{0}^{t}\left(\||\partial_{1}|^\alpha\theta\|_{\dot{H}^1}^2+\||\partial_{2}|^\beta\theta\|_{\dot{H}^1}^2\right)d\tau\right]\leq C(\theta^0)\|\theta^0\|_{H^2}^2\exp\left[ C(\alpha,\beta,\theta^0)\right] .$$
Finally we get
\begin{align}
	\|\theta(t)\|^2_{\dot{H}^2}+\int_{0}^{t}\||\partial_{1}|^{\alpha}\theta\|^2_{\dot{H}^2}d\tau+\int_{0}^{t}\||\partial_2|^{\beta}\theta\|^2_{\dot{H}^2}d\tau\leq C(\alpha,\beta,\theta^0),
\end{align}
Which prove the result.
\end{preuve}
\newpage
Now, we are ready for the proof of our first main result:

\begin{preuve1}
By \eqref{Hs}, we have, for any $t\geq 0$
\begin{equation*}\label{4.10}
	\|\theta(t)\|^2_{H^s}+\int_0^t||\partial_{1}|^{\alpha}\theta\|^2_{H^s}d\tau+\int_{0}^t\||\partial_2|^{\beta}\theta\|^2_{H^s}d\tau\leq \|\theta^0\|^2_{H^s}+ C\int_{0}^t\left(1+\|\nabla u\|_{L^\infty}+\|\nabla \theta\|_{L^\infty}\right) \| \theta\|_{H^s}^2d\tau,
\end{equation*}
The fact that $H^r(\R^2)\hookrightarrow L^\infty(\R^2)$, for all $r>1$, so we have
\begin{align}
	\nonumber	\|\nabla \theta\|_{L^\infty}\leq \|\nabla\theta\|_{H^{1+\min\{\alpha,\beta\} }}
	&\leq C\left(1+\||\partial_{1}|^{2+\min\{\alpha,\beta\}}\theta\|_{L^2}+\||\partial_{2}|^{2+\min\{\alpha,\beta\}}\theta\|_{L^2}\right)\\
	\label{4.11}	&\leq C\left(1+\||\partial_{1}|^\alpha\theta\|_{\dot{H}^2}^2+\||\partial_{2}|^\beta\theta\|_{\dot{H}^2}^2\right).
\end{align}
and 
\begin{align}
	\label{4.12}	\|\nabla u_\theta\|_{L^\infty}\leq \|\nabla\theta\|_{H^{1+\min\{\alpha,\beta\} }}
	&\leq C\left(1+\||\partial_{1}|^\alpha\theta\|_{\dot{H}^2}^2+\||\partial_{2}|^\beta\theta\|_{\dot{H}^2}^2\right).
\end{align}
Collecting \eqref{4.11} and \eqref{4.12} with \eqref{4.10} and the fact that
\begin{align*}
	C	\|\theta\|_{H^s}^2&\leq C\left(\||\partial_{1}|^s\theta\|_{L^2}^2+\||\partial_{2}|^s\theta\|_{L^2}^2\right)\\
	&\leq C\left(\||\partial_{1}|^\alpha\theta\|_{L^2}^2+\||\partial_{2}|^\beta\theta\|_{L^2}^2\right)+\frac{1}{2}\left(\||\partial_{1}|^\alpha\theta\|_{H^s}^2+\||\partial_{2}|^\beta\theta\|_{H^s}^2\right),
\end{align*}
we get
\begin{align*}
	\|\theta(t)\|^2_{H^s}+\frac{1}{2}\int_{0}^t\||\partial_{1}|^{\alpha}\theta\|^2_{H^s}d\tau+\frac{1}{2}\int_{0}^t\||\partial_2|^{\beta}\theta\|^2_{H^s}d\tau&\leq \|\theta^0\|_{H^s}^2+ C\|\theta^0\|_{L^2}^2+ C\int_{0}^t\left(\||\partial_{1}|^\alpha\theta\|_{\dot{H}^2}^2+\||\partial_{2}|^\beta\theta\|_{\dot{H}^2}^2\right) \| \theta\|_{H^s}^2d\tau.	
\end{align*}
Using Gronwall lemma we obtain
\begin{align*}
	\|\theta(t)\|^2_{H^s}&\leq C\|\theta^0\|_{H^s}^2\exp\left[C\int_{0}^t\left(\||\partial_{1}|^\alpha\theta\|_{\dot{H}^2}^2+\||\partial_{2}|^\beta\theta\|_{\dot{H}^2}^2\right) d\tau\right]\\
	&\leq C\|\theta^0\|_{H^s}^2\exp\left[C(\alpha,\beta,\theta^0)\right].
\end{align*}
Therefore, for any $t\geq 0$
\begin{align*}
	\|\theta(t)\|^2_{H^s}+\int_{0}^t\||\partial_{1}|^{\alpha}\theta\|^2_{H^s}d\tau+\int_{0}^t\||\partial_2|^{\beta}\theta\|^2_{H^s}d\tau&\leq C(\alpha,\beta,\theta^0).
\end{align*}
Finally, we get $\theta\in C_b(\R^+,H^s(\R^2))$ and $|\partial_{1}|^\alpha\theta,|\partial_{2}|^\beta\theta\in L^2(\R^+,H^s(\R^2))$.	
\end{preuve1}
	
\section{Proof of Corollary \ref{Corollary2.5}}
Let $\theta$ the global solution of \eqref{AQG}, then, $\theta\in  C(\R^+,H^s(\R^2))\cap C((0,+\infty),H^2(\R^2)),$ moreover, for any $T>0$
$$|\partial_{1}|^\alpha\theta,|\partial_{2}|^\beta\theta\in L^2([0,T],H^s(\R^2)).$$ 
Taking $t_0>0$ and considering the following system
	\begin{equation}\label{AQG2}\tag{$AQG_2$}
		\begin{cases}
			\partial_t\gamma+ u_\gamma.\nabla\gamma +|\partial_1|^{2\alpha}\gamma+ |\partial_2|^{2\beta}\gamma=0,\\
			\gamma(0)=\theta(t_0)\in H^2(\R^2).
		\end{cases}
	\end{equation}
	Therefore, by Theorem \ref{theorem1}, there exists a unique global solution of \eqref{AQG2}:
	$$\gamma\in C_b(\R^+,H^2(\R^2)),\quad |\partial_{1}|^\alpha\gamma,|\partial_{2}|^\beta\gamma\in L^2(\R^+,H^2(\R^2)).$$
	By the uniqueness of solution we have 
	$$\theta(t)=\gamma(t-t_0),\quad \forall t\geq t_0.$$
	Which implies that
	\begin{align*}
		\|\theta(t)\|_{L^\infty(\R^+,H^s)}&\leq \sup_{0\leq t\leq t_0}\|\theta(t)\|_{H^s}+\sup_{t\geq t_0}\|\theta(t)\|_{H^s}\\
		&\leq \sup_{0\leq t\leq t_0}\|\theta(t)\|_{H^s}+\sup_{t\geq0}\|\gamma(t)\|_{H^2}<+\infty
	\end{align*}
	Moreover
	\begin{align*}
		\int_{0}^{+\infty}\||\partial_{1}|^\alpha\theta\|_{H^s}^2d\tau&=\int_{0}^{t_0}\||\partial_{1}|^\alpha\theta\|_{H^s}^2d\tau+\int_{t_0}^{+\infty}\||\partial_{1}|^\alpha\theta\|_{H^s}^2d\tau\\
		&\leq \int_{0}^{t_0}\||\partial_{1}|^\alpha\theta\|_{H^s}^2d\tau+\int_{0}^{+\infty}\||\partial_{1}|^\alpha\gamma\|_{H^2}^2d\tau<+\infty,
	\end{align*}
	and by the same technique we have 
	\begin{align*}
		\int_{0}^{+\infty}\||\partial_{2}|^\beta\theta\|_{H^s}^2d\tau\leq \int_{0}^{t_0}\||\partial_{2}|^\beta\theta\|_{H^s}^2d\tau+\int_{0}^{+\infty}\||\partial_{2}|^\beta\gamma\|_{H^2}^2d\tau<+\infty.
	\end{align*}
	Finally, we get
	$$\theta\in C_b(\R^+,H^s(\R^2)),\quad |\partial_{1}|^\alpha\theta,|\partial_{2}|^\beta\theta\in L^2(\R^+,H^s(\R^2)).$$
\hfill $\blacksquare$

\end{document}